\documentclass[12pt, twoside, reqno]{amsart}
\usepackage{amsmath,amsthm,amsfonts,amssymb,eucal}

\renewcommand {\a}{ \alpha }
\renewcommand{\b}{\beta}
\newcommand{\y}{\eta}

\newcommand{\g}{\gamma}

\newcommand{\vark}{\varkappa}

\renewcommand{\d}{\delta}

\newcommand{\D}{\Delta}
\newcommand{\s}{\sigma}
\renewcommand{\l}{\lambda}

\renewcommand{\t}{\theta}

\newcommand{\p}{\partial}

\newcommand{\Om}{\Omega}

\newcommand{\R}{ \mathbb R}

\newcommand{\N}{ \mathbb N}
\newcommand{\Z}{ \mathbb Z}
\renewcommand{\H}{ \mathbb H}

\newcommand{\CF}{\mathcal F}
\newcommand{\CG}{\mathcal G}

\newcommand{\CP}{\mathcal P}

\newcommand{\CM}{\mathcal M}

\newcommand {\GF}{\mathfrak F}

\newcommand {\BA}{\mathbf A}

\newcommand {\BH}{\mathbf H}

\newcommand{\wt}{\widetilde}

\DeclareMathOperator{\dom}{Dom}

\DeclareMathOperator{\vol}{\rm{vol}}

\newtheorem{thm}{Theorem}[section]

\theoremstyle{definition}
\newtheorem{example}[thm]{Example}
\theoremstyle{remark}
\newtheorem{rem}[thm]{Remark}

\numberwithin{equation}{section}

\newcommand{\thmref}[1]{Theorem~\ref{#1}}

\newcommand{\bsymb}{\boldsymbol}

\newcommand{\loc}{\rm{loc}}
\newcommand{\hyp}{\rm{hyp}}

\newcommand{\bott}{\rm{bott}}

\newcommand{\sa}{self-adjoint}
\newcommand{\sh}{Schr\"odinger}

\begin{document}
\subjclass[2000]{35P15, 47F05} \keywords{Schr\"odinger
operator, Boundstates}

\title[Counting Schr\"odinger boundstates]{Counting Schr\"odinger boundstates: semiclassics and beyond}

\author[Rozenblum]{Grigori Rozenblum}
\address[Grigori Rozenblum]{Department of Mathematics \\
                        Chalmers University of Technology
                        and  The University of Gothenburg \\
                         S-412 96, Gothenburg,
                        Sweden}
\email{grigori@math.chalmers.se}
\author[Solomyak]{Michael Solomyak}
\address[Michel Solomyak]{Department of Mathematics\\ The Weizmann Institute of Science\\
Rehovot, 76100, Israel} \email{michail.solomyak@weizmann.ac.il}
\begin{abstract}
This is a survey of the basic results on the behavior of the
number of the eigenvalues of a {\sh} operator, lying below its
essential spectrum. We discuss both fast decaying potentials, for
which this behavior is semiclassical, and slowly decaying
potentials, for which the semiclassical rules are violated.
\end{abstract}
\maketitle

The outstanding personality of Sergey Lvovich Sobolev
determined the development of Analysis in XX century
in many aspects. One of his most influential
contributions to Mathematics is the invention of the
function spaces now named after him and the creation
of the machinery of embedding theorems for these
spaces. The ideology and the techniques based upon
these theorems enabled S.L. Sobolev and his followers
to find comprehensive and exact solutions  to many key
problems in Mathematical Physics. The paper to follow
is devoted to a survey of results in one of such
problems. This problem concerns the behavior of the
discrete part of the spectrum of a {\sh} operator with
negative potential.

\section{Introduction}\label{intro}
 The classical Weyl lemma states that the
essential spectrum of a self-adjoint operator $\BH$ in a Hilbert
space is stable under perturbations by a compact operator. This
lemma has many important generalizations. In particular, if $\BH$
is non-negative, the result survives if the perturbation is only
relatively compact with respect to $\BH$, in the sense of
quadratic forms.

{The leading and most inspiring example} in spectral theory, where
the Weyl lemma plays the key role, concerns the discrete spectrum
of a Schr\"odinger operator
\begin{equation*}%\label{0:1}
    \BH_V=-\D-V
\end{equation*}
 on $\R^d$. Here $V=V(x)$ is a real-valued  measurable function on $\R^d$ (the
 potential), which we assume to decay at infinity, in a certain
 appropriate { sense.  Then the operator can be defined via the corresponding
quadratic form, considered on the Sobolev space $H^1(\R^d)$. We}
assume for simplicity that $V\ge 0$.
 Results for general real-valued potentials can be then derived
 by using the variational principle.  In this paper we do not
 touch upon the results which take into account the interplay
 between the positive and the negative parts of the potential.

For the description of the spectrum of the operators involved  we
will use the following { notation}. Let $\s(\BH)$ and
 $E_\BH(\cdot)$ stand for {the spectrum} and the spectral measure of
 a {\sa} operator $\BH$. We call the number
\[ \bott(\BH):= \inf\{\l:\l\in\s(\BH)\}\]
{\it the bottom} of the operator $\BH$. We put
 \[N_-(\b;\BH)=\dim E_{\BH}(-\infty,\b),\qquad\b\in\R. \]
 The relation
\[N_-(\b;\BH)<\infty\]
means that the spectrum of $\BH$ on the half-line $( -\infty,\b)$
is discrete, {moreover, finite,} and $N_-(\b;\BH)$ gives the
number of the {eigenvalues} of $\BH$, counted according to their
multiplicities and lying on this half-line.

 \vskip0.2cm

 The spectrum of the free Laplacian $\BH_0=-\D$
 in $L^2(\R^d)$ is the half-line $[0,\infty)$, and by the Weyl lemma
 the negative spectrum of $\BH_V$ is discrete.
  However, this lemma gives no quantitative
 information about the negative spectrum: it does not allow one to
 find out, whether the quantity $N_-(0;\BH_V)$ is
 infinite, or finite, and in the latter case it gives no control
 of its size. It is often important to answer these
 questions. In order to make the problem more transparent, it is
 useful to insert a real positive parameter (the {\it coupling constant}),
 and to study the above questions for the family
\begin{equation}\label{0:2}
    \BH_{\a V}=-\D-\a V,\qquad \a>0.
\end{equation}
The function $N_-(0;\BH_{\a V})$ grows together with $\a$, and
this growth of the number of negative eigenvalues can be
interpreted as birth of new bound states from the edge of the
continuous spectrum {as the exterior field grows}. {At the same
time, $N(0;\BH_{\a V})=N(0;-\a^{-1}\D-V)$, so the behavior of this
quantity as $\a\to\infty$ describes simultaneously the
semiclassical behavior of the eigenvalues as the 'Planck constant'
$\a^{-\frac12}$ tends to $0.$}

Along with $N_-(0;\BH_{\a V})$, one often studies the
function $N_-(-\g;\BH_{\a V})$, where $\g>0$. If the
assumptions about $V$ guarantee discreteness of the
negative spectrum of $\BH_{\a V}$, then the latter
number is always finite. If, in addition,
$N_-(0;\BH_{\a V})=\infty$, the behavior of
$N_-(-\g;\BH_{\a V})$ as $\g\to 0+$, for $\a$ fixed,
is an important characteristics of the operator.

\vskip0.2cm

The main contents of the present paper is a survey and
a certain {detailing} of the known results on the
behavior of the function $N_-(-\g;\BH_{\a V})$, $\g\ge
0$, for the {\sh} operator \eqref{0:2} and its
generalizations -- such as {\sh} operators on
manifolds, or in domains $\Om\subset\R^d$. Note that
in the latter cases the bottom of the Laplacian is not
necessarily equal to zero. Then we discuss the
behavior of $N_-(\b;\BH_{\a V})$ for a fixed value of
$\b\le\bott(-\D)$ (we refrain from using the notation
$N_-(-\g,\BH_{\a V})$ except for the cases where
$\bott(\BH_0)=0$).

\vskip0.2cm

Our starting point is the Weyl asymptotic law, which allows one to
realize what sort of results is desirable.

If the potential $V$ is nice (say, $C_0^\infty$), then for any
$\g\ge 0$ the function $N_-(-\g;\BH_{\a V})$ exhibits the
semiclassical, or Weyl, asymptotic behavior, that is,
\begin{equation}\label{0:W}
    N_-(-\g;\BH_{\a V})\sim
    w_d\a^\frac{d}2\int_{\R^d}V^\frac{d}2dx,\qquad\a\to\infty.
\end{equation}
Here $w_d=v_d(2\pi)^{-d}$, where $v_d$ stands for the volume of
the unit ball in $\R^d$. {(The word 'semiclassical' is used in
order to indicate that the expression on the right-hand side in
\eqref{0:W} is proportional to the volume of the region in the
classical phase space $\R^{2d}$ where the classical Hamiltonian
$p^2-\a V(x)$ is negative.)} In particular, the asymptotic formula
\eqref{0:W} {hints}
 that for any potential $V\in L^1_{\loc}(\R^d)$ the
function $N_-(-\g;\BH_{\a V})$ cannot grow (in $\a$) slower that
$O(\a^\frac{d}2)$. {But can it grow faster?}

In this connection, the following questions arise in a natural
way.\vskip0.2cm

{\bf A. To describe the classes of potentials that guarantee the
estimate}
\begin{equation}\label{0:wide}
    N_-(-\g;\BH_{\a V})=O(\a^\frac{d}2), \qquad \a\to\infty.
\end{equation}

\vskip0.2cm

Another important question is this:

\noindent{\bf B. Suppose that for a given potential $V$ we have
\eqref{0:wide}. Does this imply the asymptotic formula
\eqref{0:W}?}

 \vskip0.2cm One more natural question:

\noindent{\bf C. What can be said about the eigenvalues for such
potentials that the negative spectrum of $\BH_{\a V}$ is still
discrete, but \eqref{0:wide} is violated?}\vskip0.2cm

In the paper we discuss {the present situation with} answers to
these questions. The answers heavily depend on the dimension. In
particular, the answer to the question {\bf B} is YES if $d\ge3$,
and it is NO if $d=1,2$.

We also discuss the analogues of these problems for the Laplacian
on a manifold and, more briefly, on domains $\Om\subset\R^d$ and
on the lattice $\Z^d$. Note that in all these cases the situation
is understood up to a much lesser extent, than for $\R^d$.

\vskip0.2cm
 The number $N_-(0;\BH_V)$ can be interpreted as the
borderline value, for $r=0$, of the quantity
\[ S_r(V)=\sum\limits_{\l_j(\BH_V)<0}|\l_j(\BH_V)|^r,\qquad r>0.\]
Estimating such sums is important for Physics, and
this is the main subject in the so called Lieb --
Thirring inequalities. In this paper we do not touch
upon this popular topic; see \cite{LW} for a survey
and \cite{FrLieb} for newer results.

\section{Operators on $\R^d,\ d\ge 3$}\label{eucl}
\subsection{The RLC estimate.}\label{>2} In the case considered, the answer
to the questions {\bf A, B} is given by the so called Rozenblum --
Lieb -- Cwikel estimate (the RLC estimate), named after the
mathematicians who gave the first independent proofs of the
result. In the form given below the result is due to Rozenblum
\cite{R1,R2}. Other authors, see \cite{L} and \cite{Cw}, did not
discuss the necessity of the condition on $V$.

\begin{thm}\label{RLC} Let $d\ge 3$. {Then there exists a
constant $C=C(d)$ such that for any $ V\in L^\frac{d}2(\R^d),\
V\ge0,$ and any $\g\ge 0$}
\begin{equation}\label{1:rlc}
    N_-(-\g;\BH_{\a V})\le C(d)\a^\frac{d}2\int_{\R^d}V^\frac{d}2dx,
\end{equation}
and moreover, the asymptotic formula \eqref{0:W} holds.

Conversely, suppose that $d\ge3$, { for a certain $V\ge 0$ the
operator $\BH_{\a V}$ is well defined for all $\a>0$,} and for
some $\g\ge0$ the function $N_-(-\g;\BH_{\a V})$ is
$O(\a^\frac{d}2)$ as $\a\to\infty$. Then $V\in L^\frac{d}2(\R^d)$,
and, therefore, estimate \eqref{1:rlc} and asymptotic formula
\eqref{0:W} are fulfilled for an arbitrary $\g\ge 0$.
\end{thm}

Evidently, estimate \eqref{1:rlc} for any $\g>0$ and
any $\a>0$ is a consequence of its particular case for
$\g=0$ and $\a=1$. Asymptotic formula \eqref{0:W} {is
proved first by elementary methods (Dirichlet --
Neumann bracketing)} for  potentials  $V\in
C_0^\infty(\R^d)$. It extends to the general case by a
machinery, known as 'completion of spectral
asymptotics' and presented in detail in the book
\cite{BS1}, see especially Lemma 1.19 there.

The proofs given by Rozenblum, by Lieb, and by Cwikel,
used different techniques. Rozenblum's approach was
based upon the Sobolev embedding theorem in
combination with Besicovitch-type covering theorem;
Cwikel applied harmonic analysis and theory of
interpolation of linear operators. Both these proofs
apply to much more general classes of operators than
just to the Laplacian, but only in the $\R^d$-setting.
The first proof which admits generalization to other
situations, say to operators on manifolds, is due to
Lieb, who used the semigroup theory, in the form of
path integrals.

Later several other proofs were suggested, including
the ones given by Fefferman \cite{F} and by Li and Yau
\cite{LY}. For us, the latter is especially
remarkable, since it shows in an extremely transparent
form the deep connection between the 'global Sobolev
inequality' and the RLC estimate. The techniques in
\cite{LY} uses semigroup theory in a somewhat more
direct way than in \cite{L}. Like Lieb's proof, it
admits far-reaching generalizations.

\section{The general RLC inequality}\label{gen}
\subsection{The approach by Li - Yau.}\label{gen:ly}
What we present below, is an abstract version of the
Li-Yau result. It was established in the paper
\cite{LevS}, whose authors aimed at finding the most
general setting in which the approach of \cite{LY}
applies. The classical notion of sub-Markov semigroup
is used in the formulation.

Let $(\Om,\s)$ be a measure space with sigma-finite measure. We
denote $L^q(\Om)=L^q(\Om,\s)$ and
$\|\cdot\|_q=\|\cdot\|_{L^q(\Om)}$. Suppose that a non-negative
quadratic form $Q[u]$ is defined on a dense in $L^2(\Om)$ linear
subset $\dom[Q]$. We assume that $Q$ is closed and that the
corresponding self-adjoint operator $\BA=\BA_Q$ generates a
symmetric, positivity preserving semigroup. In this situation we
say that the operator $\BA$ is a {\it sub-Markov generator}. We
also suppose that there exist an exponent $q>2$ and a positive
constant $K$, such that
\begin{equation}\label{1:mar}
    \|u\|_q^2\le KQ[u],\qquad \forall u\in\dom\,Q.
\end{equation}
\begin{thm}\label{1:ly}
Let $Q[u]$ be the quadratic form of a sub-Markov generator in
$L^2(\Om)$. Suppose that estimate \eqref{1:mar} is satisfied, with
some $q>2$. Let
\begin{equation}\label{1:6}
    0\le V\in L_p(\Om),\qquad p=(1-\frac2{q})^{-1}.
\end{equation}
Then the quadratic form
\begin{equation*}%\label{0:7}
    Q_{V}[u]:=Q[u]-\int_\Om V|u|^2d\s,\qquad u\in\dom[Q],
\end{equation*}
is bounded from below in $L^2(\Om)$ and closed. The negative
spectrum of the corresponding self-adjoint operator $\BA-V$ in
$L^2(\Om)$ is finite, and
\begin{equation}\label{0:8}
    N_-(0;\BA-  V)\le C(p)K^p\int_\Om V^pd\s,\qquad
    C(p)=e^2(\frac{p}2)^p.
\end{equation}
\end{thm}
We will call \eqref{0:8} {\it the general RLC inequality}.

 It is well known that for any $d$ the (minus) Laplacian on
$\R^d$ is a sub-Markov generator. The inequality \eqref{1:mar} is
satisfied if $d\ge 3$, with $q=\frac{2d}{d-2}$, so that in
\eqref{1:6} we have $p=\frac{d}2$. This is the so-called `global
Sobolev inequality', and the sharp value of the constant $K$ is
known, see, e.g., \cite{M}, inequality (3) in Section 2.3.3. So,
\thmref{1:ly} implies the RLC estimate \eqref{1:rlc}, with an
explicitly given constant. For the case $d=3$, which is the most
interesting for Physics, this constant is slightly greater than
the best value $C(3)=.116$ in \eqref{1:rlc}, known up to now. It
should be compared with the constant $w_3=.078$ in the asymptotic
formula \eqref{0:W}. This best value is given by Lieb's approach
which we discuss in the next subsection. It is worth mentioning
here that the sharp value of the constant $C(d)$ in \eqref{1:rlc},
even for $d=3$, is up to now unknown.

\vskip0.2cm

\subsection{The approach by Lieb.}\label{gen:lieb}
Below we present the main result of the paper \cite{RS}, where an
abstract version of Lieb's approach was elaborated.

Any non-negative self-adjoint operator $\BA$ in $L^2(\Om)$
generates a contractive semigroup $e^{-t\BA}$. We suppose that
this semigroup is $(2,\infty)$-bounded, which means that for any
$t>0$ the operator $e^{-t\BA}$ is bounded as acting from
$L^2(\Om)$ to $L^\infty(\Om)$. We write
\[ \BA\in \CP\]
if the semigroup $e^{-t\BA}$ is $(2,\infty)$-bounded and
positivity preserving.

Let $K(t;x,y)$ be the integral (Schwartz) kernel of $e^{-t\BA}$.
Then the function $K(t;x,x)$ is well-defined {on} $\R_+\times\Om$,
and it belongs to $L^\infty(\Om)$ for each $t>0$. We put
\begin{equation*}%\label{1:ker}
    M_\BA(t)=\|K(t;\cdot)\|_\infty.
\end{equation*}

The main result is a {\it parametric estimate}, see \eqref{1:1}
below: it involves an arbitrary function $G(z)$ of a certain
class, as a parameter. The class $\CG$ of admissible functions $G$
is defined as follows.

 The function
 $G$ is continuous, {convex,} non-negative, grows at infinity no faster
than a polynomial, and is such that $z^{-1}G(z)$ is integrable at
zero. With each $G\in\CG$ we associate another function,
\begin{equation}\label{1:param}
    g(\l)=\int_{\R_+}z^{-1}G(z)e^{-\frac{z}{\l}}dz,\qquad \l>0.
\end{equation}
\begin{thm}\label{1:lieb}
Suppose $\BA\in\CP$ is such that the function $M_\BA(t)$ is
integrable at infinity and is $o(t^{-a})$ at zero, with some
$a>0$. Fix a function $G\in\CG$, and define $g(\l)$ as in
\eqref{1:param}. Then
\begin{equation}\label{1:1}
    N_-(0;\BA-V)\le \frac1{g(1)}\int_{\R_+}\frac{dt}{t}\int_\Om
    M_{\BA}(t)G(tV(x))d\s,
\end{equation}
whenever the integral on the right is finite.
\end{thm}
Note that finiteness of the integral in \eqref{1:1} guarantees
that the relative bound of $V$ with respect to the quadratic form
of the operator $\BA$ is smaller than $1$, so that the operator
$\BA-V$ is well-defined via its quadratic form.

\vskip0.2cm

If $(\Om,\s)$ is $\R^d$ with the Lebesgue measure, and $\BA=-\D$,
then the semigroup $e^{-t\BA}$ is positivity preserving and
$(2,\infty)$-bounded, and
$M_{-\D}(t)=(2\pi)^{-\frac{d}2}t^{-\frac{d}2}$. Since $M(t)$ is a
pure power, the choice of $G\in\CG$ is indifferent, within the
value of the constant factor in the estimate. Indeed, by a change
of variables the estimate \eqref{1:1} reduces to the form
\begin{equation}\label{1:2}
    N_-(0;\BA-V)\le C(G)\int_{\R^d}V^\frac{d}2dx,
\end{equation}
where
\[C(G)=\frac1{g(1)(2\pi)^\frac{d}2}\int_0^\infty z^{-(\frac{d}2+1)}G(z)dz.\]
The assumptions about $G$ and the finiteness of $C(G)$
dictate the restriction $d\ge 3$. The optimal choice
of $G\in \CG$ was pointed out by Lieb \cite{L}.

\vskip0.5cm The relation between Theorems \ref{1:lieb} and
\ref{1:ly} is based upon the deep connection between the Sobolev
type inequality \eqref{1:mar} and the estimate
\begin{equation}\label{HeatEstPower}
    M_{\BA}(t)\le Ct^{-\frac{d}{2}};\qquad t\in(0,\infty)
\end{equation}
for the heat kernel corresponding to   the operator $\BA=\BA_Q$.
This connection was established by Varopoulos; see \cite{Var},
Sect.II.2 or \cite{Dav}, Theorem 2.4.2.

\begin{thm}\label{EquivHetSob} If the quadratic form $Q[u]$ generates a symmetric
sub-Markov semigroup on the measure space $(\Om,\s)$, then the
estimate \eqref{HeatEstPower} with $d>2$ is equivalent to the
inequality \eqref{1:mar} with $q=\frac{2d}{d-2}$.
% for all functions $u\in\dom(Q)$.
\end{thm}

So, the result of \thmref{1:lieb} yields the general RLC
inequality \eqref{0:8} and thus, is stronger than \thmref{1:ly}.
Indeed, in the general setting the behavior of the function
$M_\BA(t)$ is not necessarily expressed by the inequality
\eqref{HeatEstPower}, with the same exponent $d$ both as $t\to 0$
and $t\to\infty$. In many cases, one has
 \begin{gather}\label{diff.estim1}
    M_{\BA}(t)\le C_0 t^{-\frac{\d}{2}},\qquad t<1;\\ \label{diff.estim2} M_{\BA}(t)\le C_\infty
    t^{-\frac{D}{2}},\qquad  t>1,
\end{gather}
with $D\neq \d$. In \cite{Var} such estimates were studied for the
sub-Laplacian on nilpotent groups, and the numbers $\d,D$ were
called there dimensions at zero, resp., at infinity. We will use
these terms as well. One encounters a similar situation when
studying the Laplacian on manifolds, or on domains in $\R^d$.

If the estimates \eqref{diff.estim1}, \eqref{diff.estim2} are
known with $\d,D>2$, the eigenvalue estimates obtained from
\eqref{1:1}  vary essentially, depending on which dimension, $\d$
or $D$ is larger.

We formulate the corresponding results, not trying to find best
possible constants, however we include the coupling parameter
$\a$.

\begin{thm}\label{Thm.Diff.estim} Under the conditions of Theorem
\ref{1:lieb} suppose that the inequalities \eqref{diff.estim1},
\eqref{diff.estim2} are satisfied with some $\d,D>2$. Then the
following eigenvalue estimates hold:
\begin{equation}\label{d>D}
N_-(0;\BH_{\a V})\le C_0'\a^{\frac{\d}{2}}\int_\Om
V^{\frac{\d}{2}} d\s+C_\infty'\a^{\frac{D}{2}}\int_\Om
V^{\frac{D}{2}}d\s,
\end{equation}
if $\d\ge D$, and
\begin{equation}\label{d<D}
N_-(0;\BH_{\a V})\le C_0''\a^{\frac{\d}{2}} \int_{\a
V\ge1}V^{\frac{\d}{2}}d\s +C_\infty''\a^{\frac{D}{2}}\int_{\a
V<1}V^{\frac{D}{2}}d\s
\end{equation}
if $\d\le D$.
\end{thm}

\begin{rem}\label{rough}
If $\d\le D$, the inequalities \eqref{diff.estim1} and
\eqref{diff.estim2} imply \eqref{HeatEstPower} with $d=D$, and
hence
\begin{equation}\label{d<D:1}
    N_-(0;\BH_{\a V})\le \wt C\a^{\frac{D}{2}}
    \|V\|_\frac{D}2^\frac{D}2,\qquad D>2.
\end{equation}
It is often important that the assumption $\d>2$, appearing in
\thmref{Thm.Diff.estim}, here is unnecessary.
\end{rem}

We discuss applications of the estimates \eqref{d>D} and
\eqref{d<D} in Sections \ref{man}, \ref{beyond}, {and}
\ref{discr}.

\section{Operators on $\R^d,\ d\ge3$: non-semiclassical behavior of $N_-(0;\BH_{\a
V})$.}\label{nons,d>2} Suppose now that $d\ge 3$ but $V\notin
L^\frac{d}2(\R^d)$, though $V(x)$ vanishes as $|x|\to\infty$,
again in some appropriate sense. Then the negative spectrum of
$-\D-\a V$ is still discrete, but the RLC inequality becomes
useless.
 In this situation some estimates for the quantity $N_-(0;\BA-\a V)$ can be
obtained by using interpolation between the RLC inequality
\eqref{1:rlc} and a remarkable result, due to Maz'ya \cite{M},
Section 2.3.3. This result gives the necessary and sufficient
conditions on a weight function $V\ge 0$ for the Hardy-type
inequality
\[\int_{\R^d}V|u|^2dx\le C(V)\int_{\R^d}|\nabla u|^2
dx,\qquad\forall u\in C_0^\infty(\R^d)\] to be satisfied.

Here we present only a particular case of the general class of
estimates obtained by this approach. See \cite{BS3} for detail.

\begin{thm}\label{1:thnonrlc}
Let $d\ge 3$. Suppose that for some $q>\frac{d}2$ the potential
$V$ satisfies the condition
\begin{equation}\label{d1:2}
    \bsymb|V\bsymb|_q^q:=\sup\limits_{t>0}\left(t^q\int_{|x|^2V(x)>t}\frac{dx}{|x|^d}\right)<\infty.
\end{equation}
Then for any $\a>0$ the operator $-\D-\a V$ on $\R^d$ is bounded
from below, its negative spectrum is finite, and the estimate
\begin{equation}\label{1:3}
   N_-(0;\BH_{\a V})\le
    C(d,q)\a^q\bsymb|V\bsymb|_q^q
\end{equation}
is satisfied.
\end{thm}

The condition \eqref{d1:2} means that the function $|x|^2V(x)$
belongs to the so-called {\it weak $L^q$-space}, usually denoted
by $L^q_w$, with respect to the measure $|x|^{-d}dx$ on $\R^d$.
The functional $\bsymb|V\bsymb|_q$ is equivalent to the norm in
this space, but it does not meet the triangle inequality itself.
The space $L^q_w$ is non-separable, and it contains the usual
space $L_q$ with respect to the same measure. Replacing in
\eqref{1:3} the functional $\bsymb|V\bsymb|_q$ by the norm in
$L^q$ coarsens the estimate, and we come to the inequality
\begin{equation}\label{1:4}
    N_-(0;\BH_{\a V})\le
    C'(d,q)\a^q\int_{\R^d}V^q|x|^{2q-d}dx,\qquad 2q>d,
\end{equation}
which looks simpler than \eqref{1:3}. The estimate \eqref{1:4} was
established in \cite{EK} by a direct approach, generalizing the
one in \cite{R2}. However, \eqref{1:4} is knowingly not exact: it
is easy to see that the finiteness of the integral in \eqref{1:4}
implies
\begin{equation}\label{1:5}
    N_-(0;\BH_{\a V})=o(\a^q),\qquad \a\to\infty.
\end{equation}
Indeed, this is certainly the case for the potentials $V\in
C_0^\infty(\R^d)$. Such potentials are dense in $L^q$ with weight
$|x|^{2q-d}$, and the procedure of completion of spectral
asymptotics, mentioned in the paragraph next to \thmref{1:rlc},
shows that \eqref{1:5} extends to all $V$ from this space. This
nice reasoning is due to Birman (private communication). It easily
extends to the general situation, and it shows that {\it any
order-sharp estimate of order $q>\frac{d}2$ for the quantity
$N_-(0;\BH_{\a V})$ must involve some non-separable class of
potentials}.

In contrast to \eqref{1:4}, the estimate \eqref{1:3} is
order-sharp: say, for the potential $V$ which for large $|x|$ is
equal to

\begin{equation}\label{1:5a}
    V(x)=|x|^{-2}(\log|x|)^{-\frac1{q}},\qquad 2q>d,
\end{equation}
 the condition
\eqref{1:3} is satisfied, and for such potentials the asymptotics
\[N_-(0;\BH_{\a V})\sim c_q\a^q,\ c_q>0,\qquad\a\to\infty\]
is known, see \cite{BS3}.

\vskip0.2cm

The condition \eqref{d1:2} allows local singularities of $V$ at
the point $x=0$, which are stronger than those allowed by the
inclusion $V\in L^\frac{d}2(\R^d)$. The weight function $|x|^2$
and the measure $|x|^{-d}dx$ in \eqref{d1:2} can be replaced by
functions and measures {in} a rather wide class, see \cite{BS3}.
In particular, this allows one to control effects coming from
singularities of $V$ distributed on submanifolds in $\R^d$. For
example, suppose we are interested in the potentials with
singularities at the sphere $|x|=1$. Then, instead of \eqref{1:3},
one can use the estimate
\[N_-(0;\BH_{\a V})\le C\a^q\sup\limits_{t>0}\int\limits_{V(x)\left||x|-1\right|^\frac2{d}>
t}\frac{dx}{||x|-1|},\qquad 2q>d\ge3.\]
 Both this estimate and \eqref{1:3} are
particular cases of Theorem 4.1 in \cite{BS3}.

We do not think that a unified condition on the potential, which
is necessary and sufficient for $N_-(0;\BH_{\a V})=O(\a^q)$ with a
prescribed value of $q>\frac{d}2$, does exist.

\section{Operators on the semi-axis}\label{d1}
\subsection{Semiclassical behavior} In the case $d=1$ it is natural to deal with
the operators on the semi-axis $\R_+$, defined as
\begin{equation}\label{d1:op}
    \BH_{\a V}u(x)=-u''(x)-\a V(x)u(x),\qquad u(0)=0.
\end{equation}
An accurate definition can be given via the corresponding
quadratic form. The case of operators on the whole axis is easily
reduced to this one, by imposing the additional Dirichlet
condition at $x=0$ and adding up the two similar estimates for the
operators acting on the positive and the negative semi-axis. {The
term $+1$ must be included in the right-hand side of the resulting
estimate, since imposing this boundary} condition means the
passage to a subspace of codimension $1$ in $H^1(\R)$. Appearing
of the term $+1$ reflects the fact that $\l=0$ is a resonance
point for the operator $-\frac{d^2}{dx^2}$ in $L^2(\R)$. This
means that for an arbitrary non-trivial potential $V\ge 0$ at
least one eigenvalue exists for any $\a>0$. Hence, no estimate
homogeneous in $\a$ is possible.

The character of estimates for the operator \eqref{d1:op} is quite
different from the RLC inequality which governs the case $d\ge 3$.
The necessary and sufficient condition for the semiclassical order
\begin{equation}\label{d1:semi}
    N_-(0;\BH_{\a V})=O(\a^\frac12)
\end{equation}
is given by \thmref{d1:th1} below. However, this
condition hardly can be re-formulated in purely
function-theoretic terms.

 With any function $0\le V\in L^1_{\loc}(\R_+)$ we associate the
sequence {\bf $\bsymb{\y}(V)=\{\y_j(V)\},\ j\in\Z$,} where
\begin{equation}\label{d1:seq}
    \y_j(V)=2^j\int_{I_j}V(x)dx,\qquad I_j=(2^j,2^{j+1}),\ j\in\Z.
\end{equation}
It is not difficult to show that
\begin{equation*}%\label{d1:est0}
    N_-(0;\BH_{\a V})\le
    C\a^\frac12\sum\limits_{j\in\Z}\y_j^\frac12(V),
\end{equation*}
so that the condition
\begin{equation}\label{d1:half}
    \bsymb\y(V)\in\ell^\frac12
\end{equation}
is sufficient for the estimate \eqref{d1:semi}. It also guarantees
validity of the Weyl asymptotics, which in this case takes the
form
\begin{equation}\label{d1:W}
    N_-(0;\BH_{\a V})\sim
    \pi^{-1}\a^\frac12\int_{\R_+}V^\frac12dx,\qquad\a\to\infty.
\end{equation}
However, the condition \eqref{d1:half} is not necessary either for
\eqref{d1:semi}, or for \eqref{d1:W}.

In order to write the necessary and sufficient condition, let us
consider the family of eigenvalue problems on the intervals $I_j,\
j\in\Z$:
\begin{equation}\label{d1:eq}
    -\l u''(x)= V(x)u(x) \ {\rm {on}}\ I_j,\qquad
    u(2^j)=u(2^{j+1})=0.
\end{equation}
Here it is convenient for us to put the spectral parameter in the
left-hand side, then for each $j$ the eigenvalues $\l_{j,k},\
k=1,2,\ldots,$ of the problem \eqref{d1:eq} correspond to a
compact operator. Let $n_j(\l)$ stand for their counting function:
\[ n_j(\l)=\#\{k: \ \l_{j,k}>\l\},\qquad \l>0.\]
Each function $n_j(\l)$ satisfies the estimate
\begin{equation}\label{d1:est}
    \l^\frac12 n_j(\l)\le C\left(2^j\y_j(V)\right)^\frac12
    =C2^j\left(\int_{I_j}Vdx\right)^\frac12
\end{equation}
and exhibits the Weyl asymptotic behavior:
\begin{equation}\label{d1:as}
    \l^\frac12n_j(\l)\to\pi^{-1}\int_{I_j}V^\frac12dx.
\end{equation}
The estimate \eqref{d1:est} is uniform in $j$ (i.e., the constant
$C$ does not depend on $j$), but the asymptotics \eqref{d1:as} is
not. This is reflected in the fact that the potential $V$ is
involved in \eqref{d1:est} and in \eqref{d1:as} in two different
ways.

The following result was obtained in \cite{NS}.

\begin{thm}\label{d1:th1}
Let $0\le V\in L^1_{loc}(\R_+)$, and let $\BH_{\a V},\ \a>0$, be
the family of operators \eqref{d1:op}. The two conditions
\begin{equation}\label{d1:eta}
    \#\{j\in\Z:\y_j(V)>\l\}=O(\l^{-\frac12}),\qquad
    \l+\l^{-1}\to\infty,
\end{equation}
and
\begin{equation}\label{d1:We}
    \sup\limits_{\l>0}\sum_j\l^\frac12n_j(\l)<\infty
\end{equation}
are necessary and sufficient for the semiclassical order
\eqref{d1:semi} of the quantity $N_-(0;\BH_{\a V})$.
\end{thm}

The condition \eqref{d1:eta} means by definition that the sequence
$\bsymb\y(V)$ belongs to the {\it weak $\ell^\frac12$-space}
(notation $\ell^\frac12_w$). This condition is much weaker than
\eqref{d1:half}.

 The conditions
\eqref{d1:eta} and \eqref{d1:We} do not guarantee the
Weyl asymptotics \eqref{d1:W}. The necessary and
sufficient condition on $V$ for validity of this
asymptotics was also established in \cite{NS}; we do
not duplicate it here. Note only that in \cite{NS} a
series of examples was constructed of potentials $V$
for which the estimate \eqref{d1:semi} holds {but the
asymptotic formula is valid with the coefficient
different from {the one in \eqref{d1:W}.} This is
impossible in dimension $d\ge 3$.

\vskip0.2cm

\subsection{Non-semiclassical behavior of $N_-(0;\BH_{\a
V})$.}\label{2q>1} The situation here turns out to be rather
simple. The criterium for $N_-(0;\BH_{\a V})=O(\a^q)$ with a given
$q>\frac12$ can be expressed in terms of the same sequence
\eqref{d1:seq}.

\begin{thm}\label{d1:th2}
Let $0\le V\in L^1_{loc}(\R_+)$, and let $2q>1$. The condition
\begin{equation}\label{d1:etaq}
    \#\{j\in\Z:\y_j(V)>\l\}=O(\l^{-q}),\qquad
    \l+\l^{-1}\to\infty,
\end{equation}
is necessary and sufficient for $N_-(0;\BH_{\a V})=O(\a^q)$, and
the inequality
\begin{equation}\label{d1:estq}
    N_-(0;\BH_{\a V})\le C_q\a^q \sup\limits_{\l>0}\l^q\#\{j\in\Z:\y_j(V)>\l\}
\end{equation}
is satisfied.

The condition  similar to \eqref{d1:etaq}, with $o(\l^{-q})$ on
the right, is necessary and sufficient for $N_-(0;\BH_{\a
V})=o(\a^q)$.
\end{thm}
In particular, the condition \eqref{d1:etaq} with $q=1$ is
fulfilled, provided that
\[ \int_{\R_+} xV(x)dx<\infty.\]
The inequality
\[ N_-(0;\BH_{\a V})\le \a\int_{\R_+} xV(x)dx\]
is the classical Bargmann estimate, see, e.g., \cite{ReedSim4}.
So, the inequality \eqref{d1:estq} covers this result, within the
value of the constant factor. Note that under the Bargmann
condition one always has $N_-(0;\BH_{\a V})=o(\a)$. The argument
is the same as in Section \ref{nons,d>2}.

 \vskip0.2cm
The proof of \thmref{d1:th2} can be found in \cite{BS1}, where
actually more general multi-dimensional problems were analyzed,
and in \cite{BLS}. See also \cite{BL}, where the result is
presented without proof.

\thmref{d1:th2} turns out to be quite useful for the
estimation {of} $N_-(\g;\BH_{\a V})$ for such
multi-dimensional problems where an additional
`channel' can be singled out, that contributes
independently to the behavior of this function for
large values of $\a$. This happens, for instance, in
many problems on manifolds, see Section \ref{beyond}.
Another, may be the most striking example, is
connected with the Laplacian on $\R^2$. We discuss
this case  in the next section.

\section{Operators on $\R^2$}\label{d2}
\subsection{Semiclassical behavior}\label{d2:semi}
In the borderline case $d=2$ the exhaustive description of the
class of potentials such that $N_-(0;\BH_{\a V})=O(\a)$, or at
least
\begin{equation}\label{d2:1}
    N_-(-\g;\BH_{\a V})=O(\a),\qquad\forall \g>0,
\end{equation} is not
known up to present. On the technical level, this is a consequence
of the fact that the embedding theorem $H^1(\R^d)\subset
L^q(\R^d),\ q=\frac{2d}{d-2}$, fails for $d=2$ (when $q=\infty$),
or of the equivalent fact that $M_\D(t)=ct^{-1}$, see Subsection
\ref{gen:lieb},
 {so the integral in \eqref{1:1} diverges}. There are various
sufficient conditions on the potential which ensure the
order-sharp in $\a$ estimate for the function \eqref{d2:1}, but
all of them are not sharp in the function classes for $V$. Even
the most general sufficient condition of this type, known up to
now (formulated in terms of Orlicz spaces), see \cite{S-dim2}, is
not necessary. What is more, there are problems of a rather close
nature, for which the RLC-like condition $V\in L^1(\R^2)$ turns
out to be sufficient, see \cite{Lap, LapNet}. So, for $d=2$ the
situation is not well understood up to now.

Below we present a comparatively simple sufficient condition,
which was {found} in \cite{BL}. Fix a number $\vark>1$, and with
any potential $0\le V\in L^\vark_{\loc}(\R^2)$ let us associate
the sequence $\bsymb\t(V,\vark)=\{\t_j(V,\vark)\},\ j=0,1,\ldots$,
where
\begin{equation*}%\label{d2:2}
    \t_0(V,\vark)^\vark =\int_{|x|<1}V^\vark dx,
    \end{equation*}
\begin{equation*}%\label{d2:2x}
    \t_j(V,\vark)^\vark =\int\limits_{2^{j-1}<|x|<2^j}|x|^{2(\vark-1)}V^\vark
    dx,\qquad
    j\in\N.
\end{equation*}
\begin{thm}\label{d2:thm1}
{For any  fixed numbers $\vark>1$ and  $\g>0$ there exists a
constant $C(\g,\vark)>0$ such that as soon as
\begin{equation}\label{d2:2a}
    \bsymb\t(V,\vark)\in\ell^1,
\end{equation}
 the operator $\BH_{\a V}$ is bounded
from below  for any $\a>0$ , its negative spectrum is discrete,
and}
\begin{equation*}%\label{d2:3}
    N_-(-\g;\BH_{\a V})\le C(\g,\vark)\a\|\bsymb\t(V,\vark)\|_1.
\end{equation*}
\end{thm}

The constant $C(\g,\vark)$ may blow up as $\g\to 0$, and the
assumption \eqref{d2:2a} does not guarantee the semiclassical
behavior of $N_-(0;\BH_{\a V})$. It turns out that for the
analysis of this behavior one has to consider separately the
subspace $\CF$ of radial functions, $u(x)=f(|x|)$. On $\CF$ the
quadratic form of $\BH_{\a V}$ generates {a second order} ordinary
differential  operator  whose spectrum is not controlled by the
sequence \eqref{d1:2}. In order to control it and to have the
semiclassical order $N_-(0;\BH_{\a V})=O(\a)$ for the original
operator, one uses \thmref{d1:th2} with the exponent $q=1$. In the
next theorem we present the final result which can be obtained by
means of this approach; in formula \eqref{d2:4} below we express
the potential $V$ in the polar coordinates.

\begin{thm}\label{d2:thm2}
Let $V\ge 0$ satisfy the conditions of \thmref{d2:thm1}. Consider
an auxiliary `effective potential' on $\R_+$,
\begin{equation}\label{d2:4}
    F_V(t)=e^{2t}\int_{-\pi}^\pi V(e^t,\phi)d\phi,\qquad t>0.
\end{equation}
Let $\{\y_j(V)\}$ denote the sequence \eqref{d1:seq} for the
potential $F_V$. Then \eqref{d2:1} holds also for $\g=0$ if and
only if the additional condition \eqref{d1:etaq} with $q=1$ is
fulfilled.
\end{thm}

Note that {by changing \eqref{d1:etaq} to} a stronger condition,
with $o(\l^{-1})$ {on} the right, we come to a condition ensuring
the Weyl asymptotics \eqref{0:W} for $d=2$. See the papers
\cite{S-dim2} and, especially, \cite{BL} for more detail and for
discussion.

This effect (appearance of an additional differential operator in
a lower dimension, which contributes to the behavior of
$N_-(0;\BH_{\a V})$ in an independent way) we meet in several
other problems, discussed in Section \ref{beyond}. This can be
interpreted as opening of an additional channel which affects the
behavior of the system studied.

\subsection{Non-semiclassical behavior}\label{d2:nons}
It is easy to see that the condition
$\bsymb{\t}(V,\vark)\in\ell^\infty$ is sufficient for
{form-}boundedness in $H^1(\R^2)$ of the multiplication by $V$.
The next result follows from {this property} and from Theorems
\ref{d2:thm1}, \ref{d2:thm2} by interpolation.

\begin{thm}\label{d2:thm3}
1. Suppose that for some $q>1$ we have
\begin{equation*}%\label{d2:5}
    \#\{j\in\N:\t_j(V,\vark)>\l\}=O(\l^{-q}).
\end{equation*}
Then for any $\g>0$ and $\a>0$
\begin{equation}\label{d2:6}
    N_-(-\g;\BH_{\a V})\le
    1+C_{\g,q}\a^q\sup\limits_{\l>0}\l^q\#\{j\in\N:\t_j(V,\vark)>\l\}.
\end{equation}

2. Besides, suppose that the sequence $\y_j(V)$, introduced in
\thmref{d2:thm2}, satisfies the condition \eqref{d1:etaq}, with
the same value of $q$. Then $N_-(0;\BH_{\a V})=O(\a^q)$, and the
function $N_-(0;\BH_{\a V})$ is controlled by the expression as in
the left-hand side of \eqref{d2:6}, plus the additional term
\[\a^q\sup\limits_{\l>0}\l^q\#\{j\in\Z:\y_j(V,\vark)>\l\}.\]
\end{thm}

\section{Schr\"odinger operator on manifolds}\label{man}
\subsection{Preliminary remarks}
Let $\CM=\CM^d$ be a smooth Riemannian manifold of dimension
$d\ge3$, and let $dx$ stand for the volume element on $\CM$. In
this section we discuss the behavior of the function
$N_-(\b;-\D_\CM-\a V)$ where $\D_\CM$ is the Laplacian on $\CM$
(i.e., the Laplace-Beltrami operator). In order to avoid any
ambiguity, here we do not use the shortened notation $\BH_{\a V}$.
As a rule, we suppose that $\CM$ is non-compact. Otherwise, the
spectrum of $-\D_\CM$ is discrete, and it makes no sense to speak
about birth of eigenvalues of $-\D_\CM-\a V$ from the essential
spectrum of $-\D_\CM$.

For a complete Riemannian manifold $\CM$ the operator $-\D_\CM$,
defined initially on $C^\infty_0(\CM)$, is essentially
self-adjoint and generates a sub-Markov semigroup. Thus, the
results of Theorems \ref{1:ly} and \ref{1:lieb} can be applied as
soon as one has sufficient information about the embedding theorem
on $\CM$, or about estimates of the heat kernel. The global
Sobolev inequality \eqref{1:mar} with the correct order
$q=2d(d-2)^{-1}$ holds only in some special cases, and for general
manifolds, probably, the only existing approach is based upon heat
kernel estimates of the type \eqref{diff.estim1},
\eqref{diff.estim2}. Usually (though, not always)
\eqref{diff.estim1} is satisfied with $\d=d$. For example, this is
the case for the manifolds of bounded geometry, see, e.g.,
\cite{Grig}. In the discussion below we will assume that
\begin{equation}\label{man:d}
    M_{-\D_\CM}(t)\le C_0 t^{-\frac{d}2},\qquad t<1.
\end{equation}
On the other hand, $D$ in \eqref{diff.estim2} reflects the global
geometry of the manifolds, however, rather roughly, and any
relation $d>D,\ d=D$, or $d<D$ is possible. The results that
follow from such estimates are given by \thmref{Thm.Diff.estim},
where one should take $\d=d$.

An important difference from the case $\CM=\R^d$ is that now the
possibility {of a positive } $ \b_\CM:=\bott(-\D_\CM)$
 is not excluded. May be, the only general result which holds true
 for any manifold {subject to \eqref{man:d}} is the following elementary,
 but useful statement.

 \begin{thm}\label{man:th1}
{Let $\CM=\CM^d$ be  a smooth complete Riemannian manifold, $d\ge
3$.} Suppose that  the inequality \eqref{man:d} is satisfied. Then
for any $0\le V\in L^\frac{d}2(\CM)$ and for any $\b<\b_\CM$ the
following inequality holds:
\begin{equation}\label{man:2}
    N_-(\b;-\D_\CM-\a V)\le C(\CM,\b)\a^\frac{d}2\int_\CM V^\frac{d}2dx,
    \qquad\forall\a>0.
\end{equation}
Along with the estimate \eqref{man:2}, the Weyl asymptotic formula
\[ N_-(\b;-\D_\CM-\a V)\sim w_d\a^\frac{d}2\int_\CM V^\frac{d}2dx,\qquad
\a\to\infty\] is satisfied.
\end{thm}

We only outline the proof of inequality \eqref{man:2}. For any
$\b<0$, the semigroup $e^{-t(-\D_{\CM}-\b)}$ is sub-Markov
(together with $e^{t\D_\CM}$), and the function
$M_{-\D_\CM-\b}(t)=e^{\b t}M_{\D_\CM}(t)$ satisfies the same
estimate \eqref{man:2}. Besides, this function exponentially
decays as $t\to\infty$, and hence \eqref{diff.estim2} is fulfilled
with {any}  $D$. So, applying \eqref{d<D:1} {with} $D=d$ to the
semigroup generated by the operator $-\D_{\CM}-\b$, we justify
\eqref{man:2} for any $\b<0$. It extends to any values $\b<\b_\CM$
by the standard variational argument. One should {only} take into
account that for all $\b<\b_\CM$ the quadratic forms
\[ \int_\CM(|\nabla u|^2-\b|u|^2)dx\]
generate mutually equivalent metrics on the Sobolev space
$H^1(\CM)$.

\vskip0.2cm

The main issue in this type of problems is whether the estimate
\eqref{man:2} remains valid for $\b=\b_\CM$. Just such an
estimate, rather than \eqref{man:2} for $\b=0$, should be
considered as the genuine generalization of the RLC inequality
\eqref{1:rlc} to the operators on manifolds. The answer to this
question is positive only in some special cases. The Hyperbolic
Laplacian is one of these cases.

\subsection{Hyperbolic Laplacian}\label{hyp}
Let us consider the $d$-dimensional Hyperbolic space $\H^d$ for
$d\ge 3$, in the upper half-space model. This means that $\H^d$ is
realized as
 $\R^d_+:=\R^{d-1}\times\R_+$, with the metric
 \[ ds^2=z^{-2}(|dy|^2+dz^2),\qquad y\in\R^{d-1},\ z\in\R_+.\]
 The corresponding volume element is $dv_{\hyp}=z^{-d}dydz$. Recall that
 the Hyperbolic Laplacian is given by
 \[\D_{\hyp}=z^2(\D_y+\p^2_z)-(d-2)z\p_z,\]
 where $\D_y$ stands for the Euclidean Laplacian on $R^{d-1}$.
 The {bottom}  of $-\D_{\hyp}$ is the point
 $\b_0(d)=\frac{(d-1)^2}4$.

 The following result, which can be called the RLC estimate for the
 Hyperbolic Laplacian, was obtained in \cite{LevS}.
 \begin{thm}\label{rlchyp}
 Let $d\ge3$ and $0\le V\in L^\frac{d}2(\H^{d})$. Then
 \[N_-(\b_0(d);-\D_{\hyp}-\a V)\le
 C(d)\a^\frac{d}2\int_{\H^d}V^\frac{d}2dv_{\hyp}.\]
\end{thm}
For the proof, one considers the quadratic form  of $-\D_{\hyp}$
which is
\[Q[u]=\int_{\H^d}(|\nabla_{\hyp}u|^2-\b_0(d)|u|^2)dv_{\hyp}=
\int_{\R^d_+}(\nabla u|^2-\b_0(d)|u|^2)z^{2-d}dydz.\] The function
$\phi(y,z)=z^{\frac{d-1}2}$ satisfies the equation
$-\D_{\hyp}\phi=\b_0(d)\phi$. The standard substitution $u=w\phi$
reduces $Q[u]$ to the form
\[ Q[u]=\int_{\R^d_+}|\nabla w|^2zdydz.\]
For this quadratic form the lower bound is already $\b=0$. Now the
global Sobolev inequality, which allows to apply \thmref{1:ly} and
leads to the estimate in \thmref{rlchyp}, follows from \cite{M},
Corollary 2.1.6/3.

\section{Operators on manifolds: beyond \thmref{Thm.Diff.estim}.}\label{beyond}
A theory, allowing one to describe the potentials $V$
on a general manifold, which ensure the semiclassical
behavior $N_-(\b_\CM;-\D_\CM-\a V)=O(\a^\frac{d}2)$,
does not exist. The situation simplifies if one has a
more detailed information about the manifold, than
that given by the values of the exponents $\d$ and $D$
in the inequalities \eqref{diff.estim1},
\eqref{diff.estim2}. We illustrate this by several
examples. We start with the simple case of a compact
manifold.

\begin{example}\label{ex.compact}
Let $\CM$ be a compact and connected Riemannian
manifold {of dimension $d\ge3$}. Then the spectrum of
$\BA=-\D_\CM$ is discrete. The number $\l_0=0$ is a
simple eigenvalue of $-\D_\CM$, the corresponding
eigenspace $\CF$ is formed by constant functions on
$\CM$. So, we have $\b_{-\D_\CM}=0$. The estimate
\eqref{man:2} for $\b=0$ certainly fails, which
immediately follows from the analytic perturbation
theory: indeed,  it shows that for any non-trivial
$V\ge 0$ and any $\a>0$ the operator $-\D_\CM-\a V$
has at least one negative eigenvalue. On the contrary,
\eqref{man:2} with $\b=0$ would give $N_-(0;-\D_\CM-\a
V)=0$ for $\a$ sufficiently small.

It is easy to show that instead of  \eqref{man:2} we have in this
example:
\begin{equation}\label{compact}
    N_-(0;-\D_{\CM}-\a V)\le 1+C(\CM) \a^{\frac{d}{2}}\int_\CM V^{\frac{d}{2}}
    dx.
\end{equation}
\end{example}
The estimate \eqref{compact} has the same properties as the RLC
estimate for $\R^d$: it gives the correct order in $\a\to\infty$
and it involves the sharp  class of potentials for which this
order is correct.

In general, for noncompact manifolds, one or both of these
properties can  be lost and some additional reasoning must be
used.

In the  case $d>D>2$ the estimate \eqref{d>D} implies \eqref{0:8}
with $2p=d$ for any compactly supported $V$, and, similarly to the
case of a compact manifold, this result is sharp. Next, if the
support of $V$ has infinite measure and $V\in L^{\frac{d}{2}}\cap
L^{\frac{D}2}$, neither of the terms in \eqref{d>D} majorizes the
other one for a fixed $\a$, however when $\a\to\infty$, the first
term in \eqref{d>D} dominates. This indicates that {it} is
possible to relax the condition of finiteness of the expression in
\eqref{d>D} and still have the semiclassical order in large
coupling parameter. This difference in the dimensions $d,D$ may
generate an additional channel, which can contribute to the
behavior of $N_-(\b_\CM;-\D_\CM-\a V)$ in a non-trivial way.

In the next example $\CM$ is a product manifold.
\begin{example}\label{product}
Let $\CM=\CM_0\times\R^m$, where $\CM_0$ is a compact,
connected smooth manifold of dimension $d-m$. We
suppose that $d\ge 3$. Denote the points on $\CM$ as
$(x, y)$ where $x\in \CM_0$ and $y \in \R^m$; further,
$dx$, $dy$ stand for the volume element on $\R^m$ and
on $\CM_0$ respectively, then the volume element on
$\CM$ is $d\s= dxdy$. The heat kernel on $\CM$ is the
product of heat kernels on $\CM_0$ and $\R^m$, and
easy calculations show that here $\d=d,\, D=m$. { If
$m>2$ the estimate \eqref{d>D} applies, as soon as
$V\in L^{\frac{d}2}\cap L^{\frac{m}2}$, however this
condition on $V$ is  not sharp since the first term in
\eqref{d>D} majorizes the second one as $\a\to\infty$.
For $m\le2$ we simply cannot apply \eqref{d>D}. The
reasoning below demonstrates a typical way to handle
such situations.}

The Laplacian on $\CM_0$ has the  lowest eigenvalue $\l_0=0$,
{simple}, {with the corresponding eigenspace consisting of
constants.} Let $\l_1$ be the first nonzero eigenvalue on $\CM_0$.
Consider the orthogonal decomposition of the space $L^2(\CM)$,
\begin{equation}\label{cyl:0}
    L^2(\CM)=\GF\oplus\wt L^2(\CM),
\end{equation}
where $\GF$ consists of functions depending only on $y$, i.e.
$u(x,y)=v(y),\, v\in L^2(\R^m)$.
 Given a function $u\in L^2(\CM),$ its orthogonal projection onto $\GF$ is
$$v(y)=\frac1{\vol\, \CM_0}\int_{\CM_0}u(x,y)dx,$$
which implies that $\wt L^2(\CM)$ consists of functions
$\wt{u}(x,y)$ with zero integral over $\CM_0$ for almost all
$y\in\R^m$. The decomposition \eqref{cyl:0} is orthogonal also in
the metric of the Dirichlet integral,
\begin{equation*}%\label{cyl:1}
    \int_\CM|\nabla u|^2dxdy=\int_\CM|\nabla \wt
    u|^2dxdy+\int_{\R^m}|\nabla v(y)|^2dy.
\end{equation*}
Denote by $\wt H^1(\CM)$ the space of those $\wt u\in \wt
L^2(\CM)$ that belong to $H^1(\CM)$. On $\wt H^1(\CM)$ the metric,
generated by the Dirichlet integral, is equivalent to the standard
metric in $H^1(\CM)$:
\begin{equation}\label{cyl:2}
    \int_\CM|\nabla \wt u|^2dxdy\ge\frac12\left(\int_\CM|\nabla \wt u|^2dxdy+\l_1\int_\CM |\wt
u|^2dxdy\right),\qquad\forall \wt u\in \wt H^1(\CM).
\end{equation}
 We also
have
\begin{equation}\label{cyl:3}
    \int_\CM V|u|^2dx\le 2 \left(\int_\CM V|\wt
    u|^2dx+\int_{\R^m}
W(y)|v(y)|^2dy\right)
\end{equation}
where the `effective potential' $W(y)$ is given by
\begin{equation}\label{eff}
    W(y)=\int_{\CM_0}V(x,y)dy.
\end{equation}

The inequalities \eqref{cyl:2}, \eqref{cyl:3}, being combined with
the variational principle, show that
 \begin{equation}\label{cyl:3.1}
 N_-(0;\D_\CM-\a V)\le
N_-(-\l_1;\D_\CM-c\a V)+N_-(0;-\D_{\R^m}-c\a W),
\end{equation}
 where $c>0$ is some constant
depending only on the value of $\l_1$, and the second
term corresponds to the Schr\"odinger operator on
$\R^m$, with the potential $-c\a W(y)$. For the first
term in \eqref{cyl:3.1} we can use the estimate
\eqref{man:2}. The appearing of the second term in
\eqref{cyl:3.1} can be interpreted as opening of a new
channel in the system under consideration. For
estimating this term, we can use Theorem
\ref{1:thnonrlc}, \ref{d2:thm2}, or \ref{d1:th2},
depending on the dimension $m$. We would like to
emphasize that here we need just the estimates of
order $O(\a^\frac{d}2)$.  For the Laplacian on $\R^m$
such estimates are non-semiclassical.

 Moreover,
suppose that $V\in L^{\frac{d}2}$ but the effective
potential $W$ given by \eqref{eff} satisfies the
conditions of one of these theorems with some
$q>\frac{d}2$. Then it may happen that the second term
in \eqref{cyl:3.1} is stronger than the first one. In
particular, if $m\ge 3$ and the potential $W$ is like
in \eqref{1:5a}, this second term, in fact, gives the
correct asymptotic behavior of the function
$N_-(0,\BH_{\a V})$.\end{example}

Recall that for the operators on the half-line \thmref{d1:th2}
gives the necessary and sufficient condition for the behavior
$N_-(0;\BH_{\a V})=O(\a^q)$ with a prescribed value of
$q>\frac12$; this condition extends to the operators on the whole
line in an obvious way. So, for $m=1$ the above construction gives
more than for $m\ge 2$. Namely, it leads to the following result.

\begin{thm}\label{man:thm2}
The two conditions: $V\in L^\frac{d}2(\CM)$ and
\[ \#\{j\in\Z:2^j\int\limits_{2^j\le|y|\le
2^{j+1}}W(y)dy\}=O(\l^{-\frac{d}2}),\ \l+\l^{-1}\to\infty\] are
necessary and sufficient for the semiclassical behavior of the
function $N_-(0;-\D_\CM-\a V)$, where $\CM=\CM_0\times\R$ is a
$d$-dimensional cylinder, $d\ge3$.
\end{thm}

Note that the inclusion $V\in L^\frac{d}2(\CM)$ does not imply any
restrictions on the behavior of $W$. Actually, under some
additional assumptions about $W$ the function $N_-(0;\BH_{\a V})$
may have regular asymptotic behavior of order $\a^\frac{d}2$ but
with the asymptotic coefficient different from classical Weyl
formula.

{This kind of results can be easily extended to manifolds with
cylindric ends.}

\vskip0.2cm

In order to better understand the mechanism lying behind {such
two-term estimates}, let us consider the free Laplacian $-\D_\CM$
in Example~\ref{product}. Separation of variables shows that
$-\D_\CM$ is unitary equivalent to the orthogonal sum of the
operators $-\D_{\R^m}+\l_k,\ k=0,1,\ldots$ where $\l_k$ are the
eigenvalues of $-\D_{\CM_0}$; recall that $\l_0=0$. So, the
structure of the spectrum of $-\D_\CM$ on $[0,\l_1)$ is determined
by the $m$-dimensional Laplacian. This makes it clear, why the
behavior of the function $N_-(0;-\D_{\R^m}-\a W)$ in dimension
$m<d$ may affect the behavior of $N_-(\b_\CM;-\D_\CM-\a V)$ for
the Laplacian on a manifold of dimension $d$. The result of
\thmref{man:th1} shows that this effect does not appear for the
function $N_-(\b;-\D_\CM-\a V)$ with $\b<\b_\CM$.

This can be considered as manifestation of the `threshold effect'
in this type of problems.

This effect exhibits in many other problems. One of them concerns
the behavior of $N_-(0;-\D-\a V)$ on $\R^2$, discussed in Section
\ref{d2}. Note that this is the problem where the effect of
appearance of an additional channel was observed and explained for
the first time, see \cite{S-dim2} and \cite{BL}. It is worth
noting also, that in the latter problem the mechanism behind this
effect is rather latent. Indeed, unlike in Example \ref{product},
here removing the `bad' subspace of radial functions does not lead
to the shift of the spectrum of the unperturbed operator.

Another class of problems where the threshold effect has to be
taken into account, concerns various periodic operators, perturbed
by a decaying potential. In this connection, see the papers
\cite{BLS, BS4}.

\vskip0.2cm

One meets similar effects when studying the behavior of
$N_-(\b;-\D_\Om-\a V)$, where $\D_\Om$ is the Dirichlet Laplacian
in an unbounded domain $\Om\subset\R^d$. For the corresponding
heat kernel the estimate \eqref{diff.estim1} with $\d=d$ always
holds. Again, it may happen that the bottom of $-\D_\Om$ is a
point $\b_0>0$. Suppose $d\ge 3$, then for any $\b<\b_0$
\thmref{man:th1} applies. So, the problem consists in finding the
estimates and the asymptotics of $N_-(\b_0;-\D_\Om-\a V)$. The
general strategy here is the same as for manifolds, and examples
like \ref{product} can be easily constructed.

\section{{\sh} operator on a lattice}\label{discr}
The techniques based upon Theorems \ref{1:ly} and \ref{1:lieb}
 applies also to the discrete Laplacian. Below we present some
 results for the simplest case, when the underlying measure space
 $(\Om,\s)$ is $\Z^d$ with the standard counting measure, so that
 $\s(E)=\# E$ for any subset $E\subset\Z^d$. For definiteness, we
 discuss only the case $d\ge 3$. The discrete Laplacian is
\[ (\BA_d u)(x)=\sum_j(u(x+1_j)+u(x-1_j)-2u(x)), \qquad x\in \Z^d,\]
 where $1_j$ is the multi-index with all zero entries except $1$ in the
position $j$. This is a bounded operator, and its spectrum is
absolutely continuous and coincides with the segment $[0,2d]$. The
corresponding heat kernel can be found explicitly, it is bounded
as $t\to 0$ and is $O(t^{-\frac{d}{2}})$
 as $t\to\infty$, thus $\d=0$ and $D=d$. The inequality
 \eqref{d<D:1}  applies, and  we obtain the discrete
 RLC estimate,
 \begin{equation*}%\label{discr:1}
    N_-(0;\BA_d-\a V)\le
    C\a^{\frac{d}{2}}\int_{\Z^d}V^\frac{d}2d\s,\quad\forall\a>0;\qquad d\ge3.
\end{equation*}
On the contrary to the continuous case, this estimate
can not be order-sharp, since the assumption $V\in
L^{\frac{d}2}(\Z^d)$ immediately yields
\[N_-(0;\BA_d-\a V)=o(\a^{\frac{d}2}),\qquad\a\to\infty.\]
Indeed, this is certainly true for any $V$ with bounded support,
since for such $V$ the number $N_-(0;\BA_d-\a V)$ is no greater
than the number $\#\{x\in\Z^d:V(x)\neq 0\}$. The set of all such
$V$ is dense in $L^{\frac{d}2}(\Z^d)$. Therefore, the result
extends to all $V\in L^{\frac{d}2}$.

We do not know even a single example of a potential $V$ on $\Z^d$,
such that $N_-(0;\BA_d-\a V)=O(\a^{\frac{d}2})$ but $\neq
o(\a^{\frac{d}2})$.

\vskip0.2cm

One more important difference with the continuous case
is that for the discrete operators the behavior
$N_-(0;\BA_d-\a V)=O(\a^q)$ with $2q<d$ is possible;
in the continuous case it never occurs in dimensions
$d\ge 3$ and, probably, also in $d=2$. In $d=1$ the
order $O(\a^q)$ with $2q<1$ is possible,  if  one
allows  potentials which  are distributions  supported
by  a subset of  zero Lebesgue  measure.

For a  given  potential $V\ge 0$ on $\Z^d$, one cannot formally
use \eqref{d<D} with $\d=0$, since the value $\d=0$ lies outside
the set admissible by \thmref{Thm.Diff.estim}. However, by using
the variational principle and \eqref{d<D} written for the
potential $V$ restricted to the set $\{x: \a V(x)<1\}$,
 it is not difficult to show that in
this particular case  \eqref{d<D} holds for any $\a>0$ even with
$\d=0$. If we introduce the distribution function of $V$,
 \[
m(\tau)=\#\{x\in\Z^d:V(x)>\tau\},\qquad \tau>0,\] this line of
reasoning leads to the inequality
\begin{equation}\label{9:1}
    N_-(0;\BA_d-\a V)\le C\left(m(2\a^{-1})+\a^\frac{d}2\int_{\a
V(x)<1}V^\frac{d}2d\s\right)
\end{equation}

{By estimating the integral in \eqref{9:1}, we come to the
following result, which has no continuous analogue.}
\begin{thm}\label{9:thm1}
{Suppose $d\ge 3$. Then for any $\nu>2$ the estimate holds}
\begin{equation}\label{9:2}
    N_-(0;\BA_d-\a V)\le C(d,\nu)\a^{\frac{d}{\nu}}\sup\limits_{\tau>0}
    \left(\tau^{\frac{d}{\nu}}\#\{x\in\Z^d:V(x)>\tau\}\right).
\end{equation}
\end{thm}
The class of discrete potentials $V(x)$, for which the functional
on the right-hand side of \eqref{9:2} is finite, is nothing but
the cone of all positive elements in the `weak' space
$\ell_w^\frac{d}{\nu}(\Z^d)$. The assumption $\nu>2$ yields
\[\ell_w^\frac{d}{\nu}(\Z^d)\subset L^\frac{d}2(\Z^d),\]
so that the estimate \eqref{9:1} applies. The inequality
\eqref{9:2} gives a better estimate than \eqref{9:1}, and it is
possible to show that, unlike \eqref{9:1}, it is order-sharp.
\vskip0.2cm

\thmref{9:thm1} applies to the potentials { decaying no slower
than $|x|^{-\nu}$, $\nu>2$, and gives the order
$O(\a^\frac{d}{\nu})$. For potentials {decaying more slowly (but
still faster that $|x|^{-2}$),}  so that the integral in
\eqref{9:1} diverges, the following result applies.  It is the
direct analogue of \thmref{1:thnonrlc}; its proof is {also} based
upon the interpolation theory.}
\begin{thm}\label{9:thm2}
\textbf Let $d\ge 3$ and $2q>d$. Suppose the potential $V(n)\ge 0$
is such that
\[
\bsymb|V\bsymb|_q^q:=\sup\limits_{\tau>0}\left(\tau^q\int_{|x|^2V(x)>\tau}
\frac{d\s}{|x|^d}\right)<\infty.\] Then the estimate
\eqref{1:3} holds for the operator $\BH_{\a
V}=\BA_d-\a V$.
\end{thm}

In connection with Theorems \ref{9:thm1} and
\ref{9:thm2} we note that for any $\nu>0$ the
potential $V(x)=(|x|+1)^{-\nu}$ belongs to the class
$\ell^\frac{d}{\nu}_w(\Z^d)$, and the potential $V$
that for large $|x|$ behaves as
\[ V(x)=|x|^{-2}\log(|x|)^{-\frac1{q}},\qquad 2q>d,\]
meets the property $\bsymb|V\bsymb|_q<\infty$. {So, these theorems
embrace the cases of estimates of orders, respectively, smaller
and larger than $\a^{\frac{d}{2}}$. It is unclear at the moment
whether a sharp estimate of the order $\a^{\frac{d}{2}}$ is
possible. This indicates that the notion of 'semiclassical' order
is not applicable here.}}

{The above results can be extended to combinatorial Schr\"odinger
operators on arbitrary infinite graphs, as soon at the heat kernel
estimates \eqref{diff.estim1}, \eqref{diff.estim2} are known with
$\d=0$, $D>2$.}

 \section{Some unsolved problems}\label{uns}
In this concluding section we list some problems in this field,
which remain unsolved up to present. In our opinion, their
solution would be important for the further progress in the field.

\vskip0.2cm

In the first place, this is the study of the Schr\"odinger
operator on $\R^2$. Here we mean an exhaustive description of
potentials $V$ ensuring the semiclassical behavior $N_-(0;\BH_{\a
V})=O(\a)$. As it was mentioned in Section \ref{d2}, the situation
here is unclear, and many natural conjectures fail to be true.
\vskip0.2cm

The next class of problems concerns manifolds. In particular, we
believe that the class of $d$-dimensional manifolds, $d\ge3$, for
which the structure of the potentials $V$, guaranteeing the
semiclassical estimate $N_-(\b_\CM;-\D_\CM-\a V)=O(\a^\frac{d}2)$,
can be exhaustively described, can be considerably widened
compared with \thmref{man:thm2}.

\vskip0.2cm

\thmref{9:thm1} indicates that the problems for the continuous and
the discrete Schr\"odinger operators have rather different nature,
and the expected results for these two parallel classes of
operators should essentially differ. It would be useful to
understand the discrete case up to a greater extent.

\vskip0.2cm

Finally, we mention the problems of the type discussed, for the
metric graphs (quantum graphs, in other terminology), in
particular for the metric trees. The few existing results, see,
e.g., \cite{EFK}, still do not give the adequate understanding of
the effects which appear when studying the Schr\"odinger operator
on graphs.

\end{document}